\documentclass{amsart}

\usepackage{amssymb} 
\usepackage{epsf}

\newcommand{\bbC}{\mathbb{C}} 
\newcommand{\bbP}{\mathbb{P}}
\newcommand{\bbR}{\mathbb{R}}
\newcommand{\Rnull}{\bbR_{\geq 0}}
\newcommand{\PC}{\bbP_2\bbC}
\newcommand{\calT}{\mathcal{T}}

\newcommand{\lra}{\longrightarrow}
\newcommand{\eps}{\varepsilon}
\renewcommand{\phi}{\varphi}

\newtheorem{lemma}{Lemma}
\newtheorem{theorem}{Theorem}
\newtheorem*{theorem*}{Theorem}
\newtheorem{cor}[theorem]{Corollary}

\newcommand{\acknowledgement}[1]{\noindent
{\it Acknowledgement\/}---#1}

\renewcommand{\tilde}[1]{\widetilde{#1}}
\newcommand{\close}[1]{\overline{#1}}
\renewcommand{\int}[1]{{#1}^\circ}

\newcommand{\Figure}[2]{%
\begin{figure}[h]

    \epsfysize#2cm
    \hspace*{\fill}\epsffile{#1.ps}\hspace*{\fill}
\end{figure} 
}
\newcommand{\FigureC}[3]{%
\begin{figure}[h]

    \epsfysize#2cm
    \hspace*{\fill}\epsffile{#1.ps}\hspace*{\fill}
    \caption{#3}
    \label{#1}
\end{figure} 
}

\begin{document}

\title[Smoothings of CQS]{Smoothings of cyclic quotient singularities
  from a topological point of view}

\author{Ludwig Balke}
\address{Math.~Institut University Bonn\\ Beringstr.~1\\53115 Bonn\\
  GERMANY}

\email{Ludwig.Balke@math.uni-bonn.de}

\subjclass{32S30 32S55 57M50}

\begin{abstract}
For a smoothing $X_s$ of a $2$--dimensional cyclic quotient
singularity $X_0 = \bbC^2/\Gamma(n,q)$, we construct a simple handle
decomposition of $X_s$ by using a particular birational map
$X_s\lra \PC$. The manifold $X_s$ is built up from the product
of an annulus with a disk by attaching $2$--handles in a manner which
can be described by means of a plumbing graph.
\end{abstract}

\maketitle

\section{Introduction}
The deformation theory of $2$--dimensional cyclic quotient
singularities was greatly clarified by results of Koll\'ar and
Shepherd-Barron. Their more general results
(cf.~\cite{Kollar-Shepherd-Barron}) implied that the reduced
components of the base space of a semiuniversal deformation of a
cyclic quotient singularity $X$ are smooth, and that there is a
bijective correspondence to certain modifications of the singularity
$X$ --- the so called P-resolutions.

The work of J.~Christophersen \cite{Christophersen} and of J.~Stevens
\cite{Stevens} gave a complete enumeration of these P-resolutions. If
$X=\bbC^2/\Gamma (n,q)$, where $\Gamma (n,q)$ is generated by
$\left(\begin{array}{cc} \zeta &\\&\zeta^q\end{array}\right)$,
$\zeta=\exp (2\pi i/n)$, with $(n,q)=1$, then the P-resolutions of $X$
are parametrized by the following combinatorial objects. We expand
$\frac{n}{n-q}$ as a continued fraction 
\[a=[a_1,\ldots,a_s]=a_1-
                        \cfrac{1}{a_2 -\cfrac{1}{\dotsb
                           - \cfrac{1}{a_s}}}, \ a_i\geq 2.\]

A sequence $k_1,\cdots,k_s$ is called a {\em chain
representing zero} (Christophersen), if all $k_i$ are positive and the
computation of the corresponding continued fraction $[k_1,\ldots,
k_s]$ yields $0$ without ever dividing by zero. With these notions,
the result of Stevens and Christophersen can be stated as follows: The
P-resolutions of $X=X(n,q)$ are parametrized by chains
$k=(k_1,\cdots,k_s)$ representing zero which have the additional
property $k\leq a$, i.e.~$k_i\leq a_i$ for $i=1,\ldots,s$.

This result gives a nice description of the smoothing components in
the base space but it says almost nothing about the topology of the
Milnor fibre for each such component.  Our main result produces an
explicit construction for the Milnor fibre in which also chains
representing zero occur.

Consider a $1$--parameter smoothing $X \lra S$ of the cyclic quotient
singularity $X_0 = \bbC^2/\Gamma (n,q)$ such that each fibre $X_s$ can
be compactified. A particular pair of canonical coordinates of
$X_0\subseteq \bbC^e$ gives a birational map $\psi_s: X_s \to \PC$. We
show that the complement of $C_s$ --- the curve where $\psi_s$
degenerates --- has a simple topological structure for $s\neq 0$: It
is homeomorphic to the product of an annulus with a disk. The manifold
$X_s$ itself is obtained from it by attaching $2$--handles
corresponding to the components of $C_s$. The manner these
$2$--handles are attached is determined by a chain representing
zero. To be more precise:

\begin{theorem*}
For $s\neq 0$, there exists a regular neighborhood of the curve
$\close{X_s}\setminus X_s$ at infinity such that the following holds:

Let $V:= \close{X_s\setminus N}$ and $B$ a regular neighbourhood in
$V$ of $C_s\cap V$. Then $A:= \close{V\setminus B}$ is diffeomorphic
to the product of an annulus with a disk and $B$ is a disjoint union
of $2$--handles.

Furthermore, $M:= A \cap \partial V$ is an oriented
Waldhausen manifold with $r$ framed boundary components and there
exists a chain $k_1,\ldots, k_s$ representing zero such that the
plumbing graph (cf. \cite{Neumann}) of $M$ is as follows:

\Figure{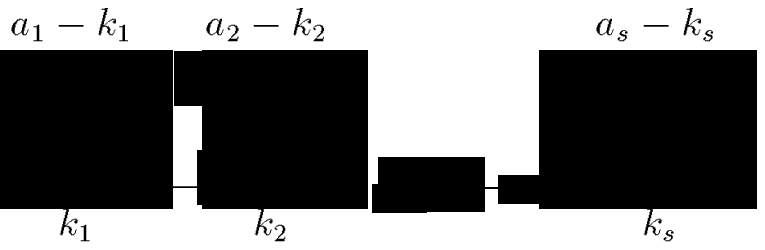}{2.5}

\noindent
We have $\frac{n}{n-q} = [a_1,\ldots, a_s], a_i \geq 2$ and
$r=\sum_{i=1}^s (a_i - k_i)$.
\end{theorem*}

This description of $X_s$ gives indeed a description of the Milnor
fibre of the smoothing due to a result of J.~Wahl (cf.~\cite{Wahl},
Thm.~2.2): Since $X_0$ is quasihomogeneous and the smoothing can be
globalized, the affine fibre $X_s$ is diffeomorphic to the interior of
the Milnor fibre of the smoothing $X\to S$ if $s$ is sufficiently
small.

As corollary of the above theorem, we can construct the homotopy type of
$X_s$ as a $CW$-complex whose $1$-skeleton is just a single $S^1$.

This paper has two parts: In the next section, we establish the
decomposition of $X_s$ into $A$ and $B$ and in the third section the
plumbing graph for $M$ is computed.

\acknowledgement{
I thank J.~Christophersen, E.~Brieskorn and O.~A.~Laudal for
stimulating discussions and helpful suggestions.
}

\section{Construction of the Decomposition}

Let $X_0 = \bbC^2/\Gamma(n,q)$, with
\[\Gamma(n,q) = 
  \left\langle \left(\begin{array}{cc} \zeta &\\&\zeta^q\end{array}\right)
                \right\rangle,\ \zeta = \exp (2\pi i/n).
\]
It is well known (cf.~\cite{Riemenschneider}) that $X_0$ can be
embedded in $\bbC^e$ where the coordinates $z_i$ are given by the
invariants $z_k = x^{i_k}y^{j_k}$ of $\Gamma (n,q)$. The numbers
$i_k,j_k$ can be computed with the help of the continued fraction
$[a_1,\ldots, a_{e-2}]= \frac{n}{n-q}$, i.e.~we have
\[a=[a_1,\ldots,a_{e-2}]=a_1-
                        \cfrac{1}{a_2 -\cfrac{1}{\dotsb
                           - \cfrac{1}{a_{e-2}}}},\ a_i\geq 2.\]
as follows. We have $(i_1,j_1) = (n,0)$ and the other exponents are
recursively defined:
\[ (i_{k-1},j_{k-1}) + (i_{k+1},j_{k+1})
                     = a_{k-1} (i_k,j_k).\]

The projection $\psi$ onto the $(z_1,z_2)$--plane is a birational map
from the projective closure $\close{X_0}$ onto the projective plane
$\PC$ which is regular on $X_0$. This follows from the fact that we
have the equations
\[ z_i z_j = p_{ij} (z_{i+1},\ldots,z_{j-1}),\ i < j-1\]
with suitable polynomials $p_{ij}$, since by an inductive argument one
sees immediately that the functions $z_3,\ldots,z_e$ are rational
functions in $z_1,z_2$.

Let us consider a $1$--parameter smoothing $X\lra S$ with $S$ the
unit-disk in $\bbC$. The fibre above $s\in S$ will be denoted by
$X_s$. Hence $X_0 = \bbC^2/\Gamma(n,q)$ and $X_s$ is smooth for $s\neq
0$. The results of Arndt (cf. \cite{Arndt,Stevens}) imply that the fibres
can be compactified and that we may assume that we have in particular
the following equations for $X$:
\begin{eqnarray}
z_i z_j = P_{ij} (s, z_{i+1},\ldots,z_{j-1}),\ i < j-1
\label{equations}
\end{eqnarray}

with a suitable polynomial $P_{ij}$ with $P_{ij}(0,\cdot) = p_{ij} (\cdot)$.
Therefore the projection $\psi_s : \close{X_s} \lra \PC$ onto the
$(z_1,z_2)$--plane is a birational map which is regular on $X_s$.

For a regular neighbourhood $N$ of the curve at infinity
$\close{X_s}\setminus X_s$ let 
\[ V := \close{\close{X_s}\setminus N}.\]  
Furthermore let $C_s$ denote the curve where $\psi_s$ degenerates. For
a suitable choice of $V$ and a suitable choice of a regular
neighbourhood $B$ of $C_s\cap V$ we will show: The components of $B$
are $2$-handles, i.e.~diffeomorphic to $D^2\times D^2$, with $D^2$ the
closed unit-disk in $\bbC$, and $A:=\close{V\setminus B}$ is a union
of a $0$-handle with a $1$-handle, i.e.~diffeomorphic to $R\times D^2$
with $R$ an annulus.

This simple handle decomposition is the reason why we have chosen this
particular birational map from $X_0$ onto $\bbC^2$ and not a generic
one.

In this section we prove the existence of the above handle
decomposition. In the next one, we investigate more precisely how the
2-handles are attached to the $4$-manifold $A$.

Since $\psi_s$ is birational, the image of $C_s$ consists
of finitely many points $p_i, i=1,\ldots, m$. The functions $z_k$
considered as rational functions in $z_1,z_2$ have poles only for
$z_1=0$. Hence we have $p_i = (0, t_i)\in \bbC^2$, where we denote the
coordinates of $\bbC^2$ by $z_1,z_2$.

The functions $z_1, z_e$ play a particular role, due to the following
lemma:

\begin{lemma}
The projection $\close{X_s}\lra \PC$ to the $(z_1, z_e)$--plane is a
ramified covering. 
\label{lemma_ramified}
\end{lemma}

\begin{proof}
The statement is easy to see for $X_0$. Since $X_s$ is a deformation
of $X_0$ it also holds for $X_s$. 
\end{proof}

In the sequel, we always assume $s\neq 0$ and hence $X_s$ smooth.
By a suitable sequence of blowing ups in the points $p_i$, we can
eliminate the points of indeterminacy of $\psi_s^{-1}$. We choose a minimal
such sequence and obtain a modification $\pi : \tilde{X_s}\lra \PC$
together with a map $\phi: \tilde{X_s} \to \close{X_s}$.
Since $\pi$ is minimal and $\tilde{X_s}$ is smooth, the restriction of
$\phi$ onto the preimage of $X_s$ is biholomorphic.

In order to define regular neighbourhoods of curves, we use the
concept of a rug function, defined by Durfee (\cite{Durfee}): Let $M$
be a semialgebraic set in $\bbR^n$ and $X\subset M$ a compact proper
semialgebraic subset of $M$. A {\em rug function} for $X$ in $M$ is a
proper semialgebraic function $\alpha : M \lra \Rnull$ with
$\alpha^{-1} (0) = X$. Durfee has shown that for $\eps$ sufficiently
small the neighbourhood $\alpha^{-1} ([0,\eps])$ of $X$ is unique up
to isotopy and independent from the particualr choice of a rug
function.

On $\tilde{X_s}$, we consider several semialgebraic functions.
First, define
\[\eta : \tilde{X_s}\lra \Rnull, p \mapsto \left\{
          \begin{array}{ll}
            \frac{1}{|z_1(p)|},& |z_1 (p)| \geq 1\\
            |z_1(p)|, & |z_1 (p)| \leq 1
          \end{array}\right.
\]

Then define $\rho_0 (p) := \min (\eta (p), \frac{1}{|z_2 (p)|}).$ By
construction, $\rho_0$ is a rug-function for the compact curve
\[ D_0 = \pi^{-1}(\PC\setminus \bbC^2) \cup V(z_1).\]
(Here, as usual, $V(f)$ denots the zero-locus of a function.)
Therefore the results of Durfee imply that 
\[ N_0 := \rho_0^{-1} ([0, \eps_0])\]
is a regular neighbourhood of $D_0$ if $\eps_0$ is sufficiently small.

On the other hand, we consider the compact curve
\[D_\infty := \tilde{X_s}\setminus \phi^{-1} (X_s) \subseteq D_0.\]
Due to Lemma~\ref{lemma_ramified}, a regular neighbourhood $N$ of
$D_\infty$ can be obtained by the rug function
\[\rho : \tilde{X_s}\lra \Rnull, p\mapsto \min (\frac{1}{|z_1(p)|},
                                             \frac{1}{|z_e(p)|}).
\]
For sufficently small $\eps > 0$, the set
\[ N:= \rho^{-1} ([0,\eps])\]
is a regular neighbourhood of $D_\infty$. Since $X_s$ is a deformation
of $X_0$ and by the uniquenes of a regular neighbourhood, we have
proven
\begin{lemma}
If $N_0'$ is a regular neighbourhood of $D_0$ then $\close{\tilde{X_s}\setminus
{N_0'}}$ is diffeomorphic to $R\times D^2$.

If $N'$ is a regular neighbourhood of $D_\infty$ then $\partial N'$ is
diffeomorphic to the link of the singularity $X_0$, i.e.~diffeomorphic
to the lens space $L(n,q)$.
\end{lemma}

Let $E$ be the union of all components of $D_0$ which are not
components of $D_\infty$, and let $\tilde{B}$ be  a regular
neighbourhood of $E$. Using the above lemma we replace $N$ by an
equivalent regular neighbourhood as follows. Define the semialgebraic
function $\alpha$ on $\tilde{X_s}$ by
\[ \alpha (p) = \left\{\begin{array}{ll}
      \rho_0 (p),&\mbox{if}\ p\notin \tilde{B}\\
      \max (\rho_0 (p),1/|z_e (p)|),&\mbox{if}\ p\in\tilde{B}
                     \end{array}\right.
\]
and redefine $N:= \alpha^{-1} ([0,\eps])$ for a sufficiently small
$\eps > 0$. Furthermore, we may assume that $\tilde{B}\setminus N =
\rho_0^{-1}([0,\eps])\setminus N$.
Then $N\cup \tilde{B}$ is a regular
neighbourhood of $D_0$. We define
\begin{eqnarray*} 
V & := & \close{\phi (\tilde{X_s}\setminus (N\cup \tilde{B}))}\\
B & := & V\cap \tilde{B}.\\
A & := & \close{V\setminus B}
\end{eqnarray*}

\FigureC{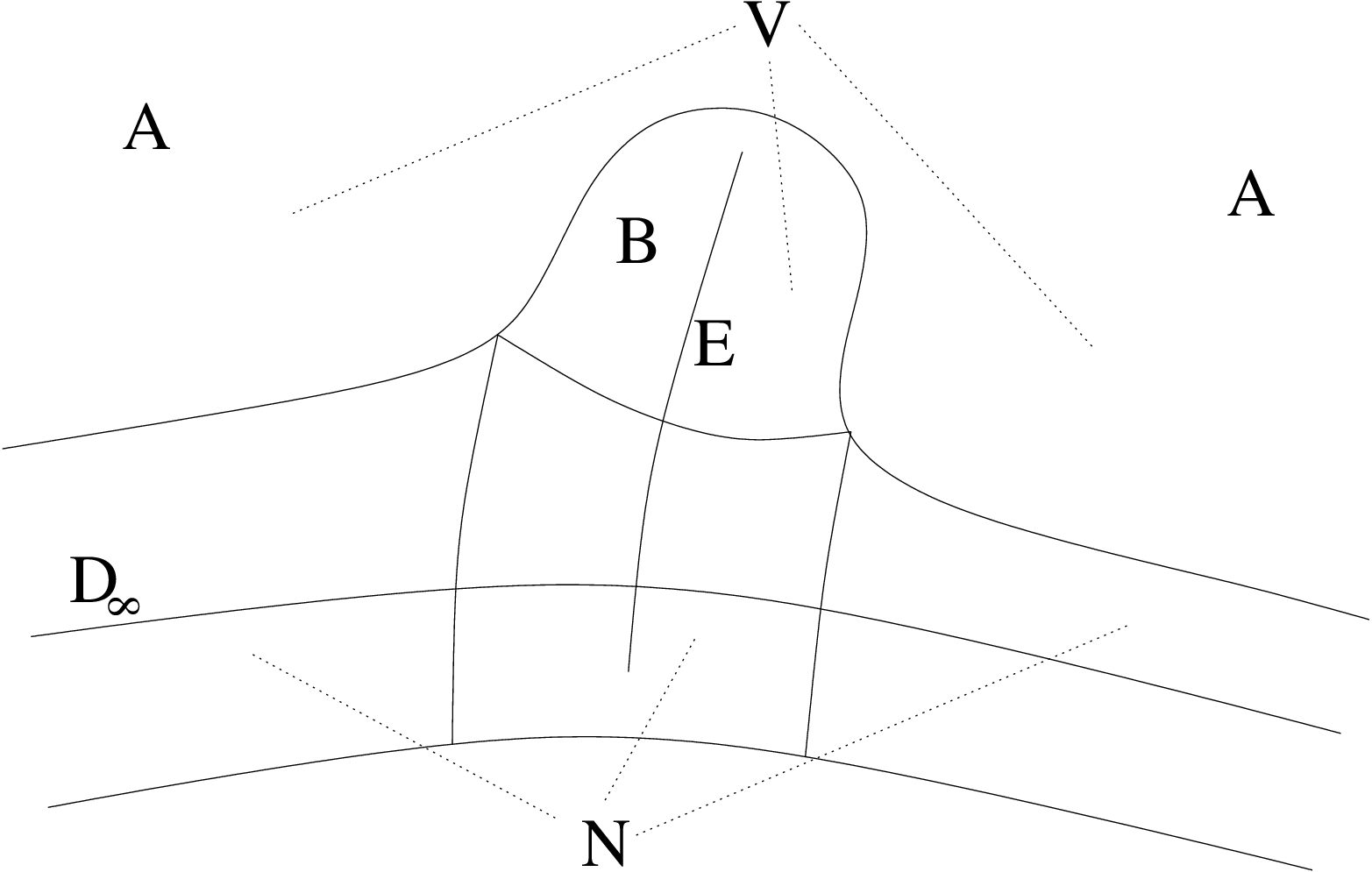}{4}{Neighbourhood of a component of $E$}

With these definitions, we want to show

\begin{lemma}
The $4$--manifold $V$ is built up from $A$, which is diffeomorphic to
$R\times D^2$, by attaching $2$-handles.

The $3$--manifold $M:= \close{\partial M \setminus B}$ is a Waldhausen
manifold.
\end{lemma}

\begin{proof}
By the above lemma, $A$ is diffeomorphic to $R\times D^2$. 
Since $\phi|_{\phi^{-1}(X_s)}$ is biholomorphic, each component of $E$
has a nontrivial intersection with $D_\infty$. Therefore it suffices
to show that for each component $F$ of $E$ the intersection $F\cap
D_\infty$ consists of a single point.

Let $B_F$ be the component of $\tilde{B}$ which contains $F$. Let $K:=
\close{F\setminus(N\cap F)}$ and $U:= \close{B_F\setminus (N\cap
  B_F)}$. Hence, $U$ is a regular neighbourhood of $K$, in particular
the fibration $\kappa : U\setminus K \lra ]0,\eps]$ is
trivialisable. Due to our assumptions on $\tilde{B}$, this implies
that the map $z_1$ defines on $U\setminus K$ a locally trivial fibration.
Therefore it is sufficient to show that a fibre of this fibration has
only one boundary component.
For this, consider the affine line $\ell_c =
\{(z_1,z_2)\in\bbC^2|\ z_1 = c\}$ with $|c|\leq \eps$. For $c\neq
0$, the function $z_e$ is holomorphic on $\ell_c$. Since $K$ is
connected, $\ell_C \cap U$ is also connected.
If the boundary of
$\ell_c\cap U$ had more than one component then we had a compact set
in $\ell_c$ such that $|z_e|$ would assume a maximum in the interior
of this compact set. But $z_e$ is not constant on $\ell_c$, therefore
$\ell_c \cap U$ has only one boundary component. This implies that
$\partial \kappa^{-1} (\eps)$ has only one component, too.

We have seen that $M$ is a particular part of a regular neighbourhood
of the curve $D_\infty$. By constructing such a regular neighbourhood
via plumbing it follows that $M$ can be decomposed into pieces which
are $S^1$-bundles over compact topological surfaces with
boundary. Hence $M$ is a Waldhausen manifold.
\end{proof}

\section{Plumbing graph of $M$}

We want to describe the manifold $M$ together with a natural framing of
its boundary components by means of a plumbing graph as introduced by
W.~Neumann. Let us shortly recall the basic properties of plumbing
graphs (for details cf.~\cite{Neumann}).
\label{section_Construction}

An oriented Waldhausen manifold is given by an oriented 3-manifold
$M$, a disjoint union $\calT$ of embedded tori and on the closure of
each component of $M\setminus \calT$ the structure of an $S^1$-bundle
over a compact surface. In particular, each boundary component is a
torus. A {\em framing} $\tau$ for $M$ is a homeomorphism defined on
the boundary such that for a component $T$ of $\partial M$ the
restriction $\tau|_T$ maps $T$ onto $S^1\times S^1$ such that the
projection on the second factor coincides with the $S^1$-bundle
structure on $T$. Up to isotopy a framing is fixed by giving the
homology classes $m:=[(\tau^{-1})_*(\{1\}\times S^1)]\in H_1 (T)$
and $f:=[(\tau^{-1})_* (S^1\times \{1\})]$. Here, $f$ is a fibre of the
$S^1$-fibration of $T$ and $m$ is a section. Moreover $(f,m)$ is an
basis of $H_1 (T)$ which gives the induced orientation on $T$.

When a framing $\tau$ is given, the manifold $M_0$, the {\em canonical
closure of $M$}, can be obtained from $M$ by gluing solid tori to the
boundary components according to the framing. The section $m$ of a
boundary component then becomes a meridian of the solid torus and the
fibre $f$ a longitude.

For our discussion it is sufficient to consider only such oriented
Waldhausen manifolds whose $S^1$-bundles are bundles over $S^2$ with a
finite number of disjoint disks removed. For such Waldhausen manifolds
the plumbing graph gives a good graphical representation of the
topological structure.

If $M$ is the manifold $S^1\times F$ with $F$ of genus $0$ with $r$
boundary components $C_1,\ldots, C_r$ together with a framing $\tau$,
then the plumbing graph $\Gamma (M,\tau)$ looks as follows:

\Figure{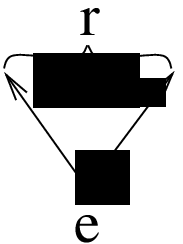}{1.7}

The arrows correspond to the boundary components $T_i = S^1\times C_i$
and we keep a bijection between the arrows and these boundary
components fixed. The weight $e$ of the vertex is the Euler number of
the $S^1$-bundle $M_0$ over $S^2$.

Let $M_1, M_2$ be oriented framed Waldhausen manifolds with plumbing
graph $\Gamma (M_1)$ and $\Gamma (M_2)$ and assume two boundary
components $T_i$ of $M_i$ given. (The case $M_1 = M_2$ is possible. In
this case we must have $T_1 \neq T_2$.) We can form a new oriented
framed Waldhausen manifold $M$ by gluing together $M_1$ and $M_2$
by identifying $T_1$ with $T_2$ in such a way that the fibre of $T_1$
becomes the section of $T_2$ and vice versa. A plumbing graph for $M$
is obtained by removing the arrows corresponding to $T_1$ and $T_2$
and joining the vertices incident with these arrows by an edge.

If $\Gamma$ is a plumbing graph for the framed oriented Waldhausen
manifold $M$, let $\Gamma_0$ be the graph which is constructed from
$\Gamma$ by removing all arrows. The resulting graph $\Gamma_0$ is
a plumbing graph for the canoncial closure $M_0$.

There are several operations on the framed Waldhausen structure which
do not change the framing or the underlying manifold but only the
decomposition of $M$ into $S^1$-bundles. W.~Neuman has given a complete
list of such operations and unique normal forms for plumbing graphs in
a more general setting than is needed here. For our purpose the
following two operations and their inverse operations are of
particalur importance:

\Figure{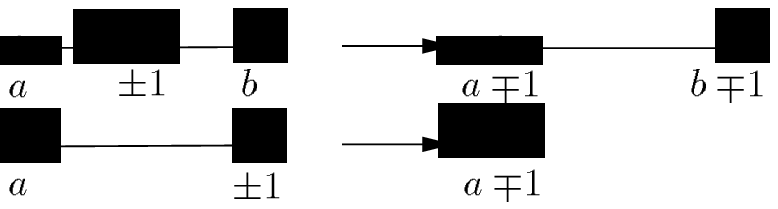}{1.7}

For obvious reasons this kind of modification is called ``blowing
down'' and its inverse ``blowing up''. This picture has to be
understood as follows. It shows the part of the plumbing graph which
is modified. All edges incident with the vertex with weight $\pm 1$
are shown, in particular there is no arrow incident with this
vertex. The other vertices may be incident with other edges or arrows
which are not shown in the picture.

For a plumbing graph $\Gamma$ let $-\Gamma$ be the graph which is
obtained by multiplying each weight of a vertex by $-1$. If $\Gamma$
has no cycles and is a plumbing graph for the oriented Waldhausen
manifold $M$ then $-\Gamma$ is a plumbing graph for the manifold $-M$
with reversed orientation. The framing of $-M$ is induced by the
framing of $M$ by replacing for each boundary component $T$ the basis
$(f,m)$ of $H_1 (T)$ by $(-f,m)$ or $(f,-m)$. After having made a
choice for one boundary there is no choice for the other components if
$\Gamma$ is connected.

The results of Neumann and the well known resolution of cyclic
quotient singularities implies the following statement.

\begin{lemma}
For the lens space $L(n,q)$ there exists up to isomoprphism exactly
one plumbing graph whose weights are all $\geq 2$, namely

\Figure{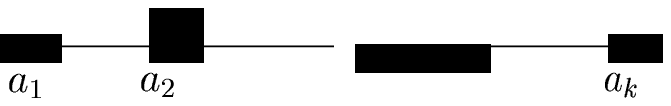}{0.7}

The numbers $a_1,\ldots,a_k$ are determined by the formula
$[a_1,\ldots,a_k] = \frac{n}{n-q}$.
\end{lemma}

\begin{proof}
First, one easily checks that $-L(n,q)$ is oriented diffeomorphic to
$L (n,n-q)$. The well known resolution of $\bbC^2/\Gamma(n,n-q)$ is a
chain of rational curves with selfintersections $-a_1,\ldots, -a_k$
where $\frac{n}{n-q} = [a_1,\ldots,a_k]$. By applying the above remark
about $-\Gamma$, we obtain a plumbing graph as desired. The uniqueness
follows from the results in \cite{Neumann}.
\end{proof}

Let us now construct a plumbing graph for the $3$-manifold $M$ defined
in the previous section. We have obtained
$\tilde{X_s}$ by successive blowing ups in the points $p_i = (0,t_i),
i=1,\ldots,m$. Let $E_i$ be the preimage of $p_i$ with respect to this
modification. Choose $\eps$ in the definition of $N$ and
$\delta > 0$ small enough such that 
\[ H_i := \{p\in\tilde{X}_s|\ |z_2(p) - t_i| \leq \delta, p\in N\cup
\tilde{B}\}\] are disjoint regular neighbourhoods of the curves
$E_i$. In particular, the boundary $\partial H_i$ of $H_i$ is
homeomorphic to $S^3$.  Let

\begin{eqnarray*}
U_i &:=& \close{H_i \cap \tilde{B}},\\
V_i &:=& \close{H_i\setminus U_i},\\
M_i &:=& M\cap V_i,\\
M_0 &:=& \close{M\setminus\bigcup_{i=1}^m M_i}
\end{eqnarray*}

Since $M$ is built up from the pieces $M_i$, $=0,\ldots,m$, we first
describe plumbing graphs for the manifolds $M_i$ with $i>0$.

$M_i$ has several boundary components. One component is the
intersection of $M_i$ with $M_0$ which we call the {\em outer}
boundary component of $M_i$. All other boundary components of $M_i$
are called {\em inner} boundary components.

Let $T$ be an inner boundary component of $M_i$, $B_T$ the
component of $U_i$ which contains $T$ and $e_T$ the component of $E$
which is contained in $B_T$. The torus $T$ 
bounds two solid tori: On the one hand $S_+ (T)$ the closure of $B_T
\cap \partial H_i$ and on the other hand $S_X (T)$ the closure of
$B_T\cap \int{H_i}$. 

\FigureC{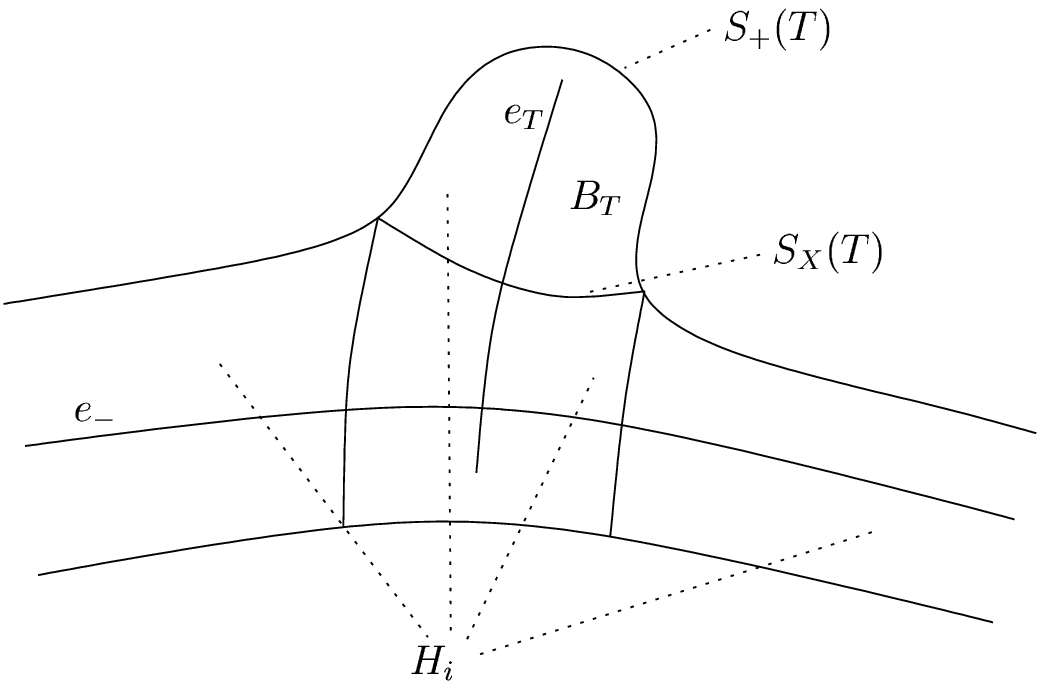}{5}{Neighbourhood of $e_T$}

We first observe:
\begin{lemma}
The curve $e_T$ has selfintersection $-1$.
\label{lemma_et}
\end{lemma}

\begin{proof}
Choose $\mu$ minimal such that $z_i$ is constant on $e_T$ for each
$i\leq \mu$. Obviously, we have $\mu \geq 2$, and according to
Lemma~\ref{lemma_ramified} we have $\mu < e$. The form of the equations
(\ref{equations}) for $X_s$ imply that the projection onto the $(z_\mu,
z_{\mu+1})$--plane is birational. Henceforth, the restriction of
$z_{\mu +1}$ to $e_T$ yields an isomorphism onto the projective line.

In the first blow up which produces $e_T$, the functions $z_\mu$ and
$z_{\mu +1}$ therefore give coordinates in a neighbourhood of
$e_T$. For that reason, there exists no point of indeterminacy of
$\psi_s^{-1}$ on $e_T$ and no further blowing up occurs in points of
$e_T$. This proves: $e_T \cdot e_T = -1$.
\end{proof}

We have two cases: There exists a compact component $e_-$ of
$D_\infty$ which intersects $e_T$ or $e_T$ intersects the strict
transform of the curve $V(z_1)$ in $\tilde{X_s}$.

Let us assume the first case. A meridian of $S_X (T)$ is a local
section of the normal bundle over $e_-$. Due to Lemma~\ref{lemma_et},
the curve $e_T$ can be blown down. After blowing down, a
meridian of $S_+ (T)$ becomes a local section in the normal bundle of
$e_-$. Together with a fibre of the bundle over $e_-$, we thus have
two framings $\tau_+$ respectively $\tau_X$ of $T$ which differ by a
fibre of the bundle.

In $B_T$, there exists exactly one compact component $e_i^0$ of $D_\infty$
which intersects the strict transform of $V(z_1)$. The normal bundle
of $e_i^0$ has this strict transform as fibre. The outer component
$T_i^0$ of $M_i$ can thus be framed in such a way that the inverse fibre of
$M_0$ is a section of $T_i^0$.

By describing $M_i$ as part of the boundary of a regular neighbourhood
of the compact curve $E_i$, we obtain in the first case a plumbing graph
$\Gamma_i$

\Figure{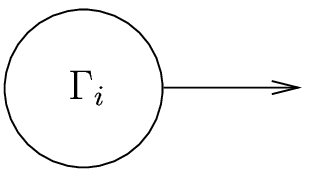}{1.7}

The arrow corresponds to $T_i^0$. Since $M_i$ as part of the boundary
of $V$ is oriented opposite to the orientation as part of the boundary
of $H_i$, the weights of the vertices of $\Gamma_i$ are all positive
numbers. For each vertex the sum of the weights and the number of
arrows incident with this vertex is $\geq 2$ since the resolution
$\pi$ is minimal.

Let $\Gamma_i'$ be the graph obtained from $\Gamma_i$ be removing all
arrows with the exception of the distinguished arrow corresponding to
$T_i^0$ and by adding $1$ for each removed arrow to the weight of
the incident vertex. Then the graph 

\Figure{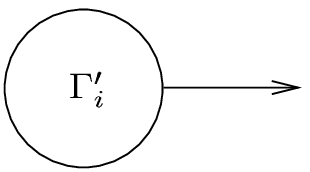}{1.7}

is a plumbing graph for the manifold $\partial V \cap H_i$.

Let us now consider the second case: There exists no compact component
of $E_i \cap D_\infty$, i.e.~one blowing up is sufficient for
eliminating the indeterminacy of $z_e$ in $p_i$. Hence $E_i$ consists
of one single component $e$ with self intersection $-1$. The manifold
$M_i$ has exactly one inner boundary component, which is parallel to
the outer component. Our arguments in the first case were local in
nature, hence we can make analogous consideration if we take as $e_-$
the strict transform of $V(z_1)$. We therefore obtain framings
$\tau_+$ and $\tau_X$ on $T$ which have as fibre the inverse fibre of
$M_0$, and the difference of the respective sections is the inverse
fibre of $M_0$.

Summarizing these constructions, we obtain a plumbing graph $\Gamma$
of $M$ as follows (possibly after a reindexing of the $M_i$):

\Figure{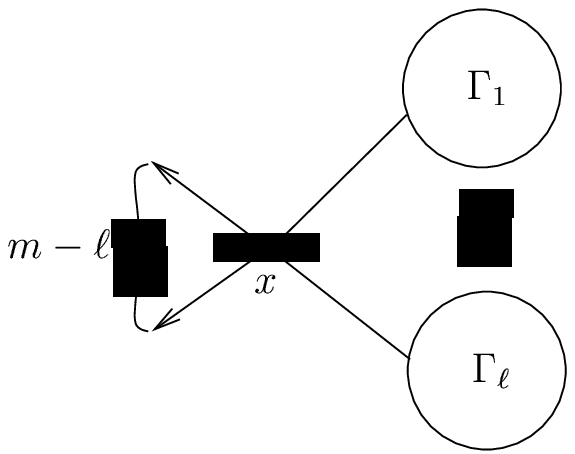}{3.8}

The graph $\Gamma_0$ which we obtain from $\Gamma$ by deleting all
arrows is then a plumbing graph for $-\partial N$ and therefore a
plumbing graph for the boundary of $R\times D^2$, i.e.~a trivial $S^1$
bundle over $S^2$. By succesive blowing down $1$-vertices we obtain
from $\Gamma_0$ the graph

\Figure{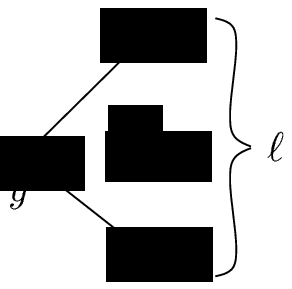}{2.5}

Therefore, we have $y\leq x$ and $y=\ell$. On the other hand, the
graph in Figure~\ref{topM3}
is a plumbing graph for $\partial V$, i.e.~for a lens space
$L(n,q)$. Since all vertices with the possible exception of the
distinguished vertex have weights $\geq 2$, we must have $\ell \leq
2$. After possibly blowing down $1$-vertices in the case $x+m-\ell =
1$ --- which is only possible for $x=y$ and $m=1$ ---
we obtain a graph $\Gamma$ for $M$ as in Fugure~\ref{topchain2}.

\FigureC{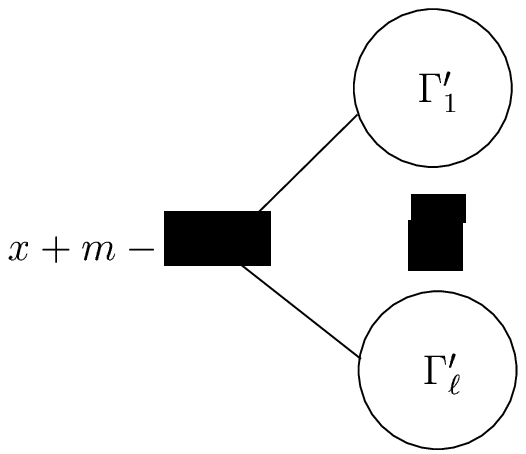}{3.5}{}

\FigureC{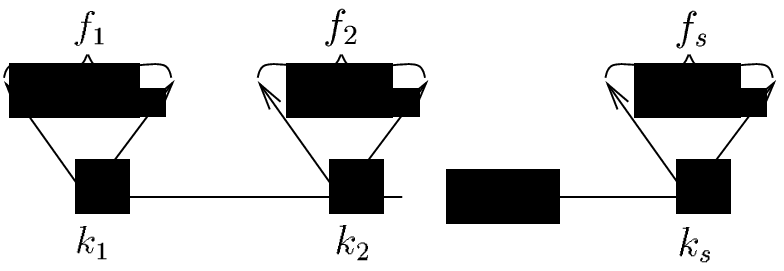}{2.5}{}

This graph has the following properties: We have $k_i + f_i \geq 2$
for $i=1,\ldots, s$. The graph $\Gamma_0$ can be blown down to a
single vertex with weight $0$. The graph $\Gamma'$ which we obtain by
removing all arrows and replacing the weight $k_i$ by $k_i+f_i$ is a
plumbing graph for the lens space $L(n,q)$. Therefore we must have
$[a_1,\ldots,a_s] = \frac{n}{n-q}$ with $a_i = k_i + f_i$, hence in
particular $s= e-2$. Since $\Gamma_0$ can be blown down to a single
vertex with weight $0$, the sequence $k_1,\ldots,k_{e-2}$ is a chain
representing zero.

Hence we can summarize our result in the following theorem:

\begin{theorem}
Let $X \lra S$ be a $1$-parameter smoothing of the cyclic quotient
singularity $X_0 = \bbC^2/\Gamma (n,q)$ such that each fibre $X_s$ can
be compactified. Let $s\in S\setminus\{0\}$.
For a suitable choice of a birational regular map $\psi_s :
\close{X_s} \lra \PC$ and a regular neighbourhood $N$ of the curve at
infinity $\close{X_s}\setminus X_s$ the following holds.

Let $V:= \close{X_t\setminus N}$ and $B$ a regular neighbourhood in
$V$ of $C_s$, the curve where $\psi_s$ degenerates. Then $A:=
\close{V\setminus B}$ is diffeomorphic to the product of an annulus
with a disk and $B$ is a disjoint union of $2$--handles.

Furthermore, $M:= A \cap \partial V$ is an oriented framed
Waldhausen manifold with $r$ boundary components whose
plumbing graph is as follows:

\Figure{topIthm}{2.5}

The sequence $k_1,\ldots,k_{e-2}$ is a chain representing zero and
$r=\sum_{i=1}^{e-2} (a_i - k_i)$.
\end{theorem}

With the notation above, we obtain immediately the following
description of the homotopy type of $V$:

\begin{cor}
The $4$-manifold $V$ is homotopy equivalent to a CW-complex with
$1$-skeleton homeomorphic to $S^1$ and $r$ cells of dimension 2.

In particular, the second Betti number of $V$ is $b_2 (V) = r - 1$ and
the fundamental group of $V$ is a finite cyclic group.
\end{cor}

\end{document}